\documentclass[12pt,reqno,oneside]{amsart}
\usepackage{amsfonts}
\usepackage{amsmath, amssymb}
\usepackage{hyperref}

\numberwithin{equation}{section} 
\newtheorem{theorem}{Theorem}[section]
\newtheorem{corollary}[theorem]{Corollary}
\newtheorem{lemma}[theorem]{Lemma}
\theoremstyle{definition}

\theoremstyle{remark}
\newtheorem{remark}[theorem]{Remark}
\numberwithin{equation}{section}

\input{tcilatex}

\begin{document}
\title[Fekete-Szeg\"{o} problem for bi-univalent functions]{Fekete-Szeg\"{o}
problem for certain classes of Ma-Minda bi-univalent functions}
\author{H. Orhan, N. Magesh and V.K.Balaji}
\address{Department of Mathematics \\
Faculty of Science, Ataturk University \\
25240 Erzurum, Turkey.\\
\texttt{e-mail:} $orhanhalit607@gmail.com$}
\address{Post-Graduate and Research Department of Mathematics,\\
Government Arts College for Men,\\
Krishnagiri 635001, Tamilnadu, India\\
\texttt{e-mail:} $nmagi\_2000@yahoo.co.in$}
\address{ Department of Mathematics, L.N. Govt College, \\
Ponneri, Chennai, Tamilnadu, India.\\
\texttt{e-mail:} $balajilsp@yahoo.co.in$}

\begin{abstract}
In the present work, we propose to investigate the Fekete-Szeg\"{o}
inequalities certain classes of analytic and bi-univalent functions defined
by subordination. The results in the bounds of the third coefficient which
improve many known results concerning different classes of bi-univalent
functions. Some interesting applications of the results presented here are
also discussed.\ \newline
\end{abstract}

\maketitle

%%%=====================

\section{Introduction}

Let $\mathcal{A}$ denote the class of functions of the form 
\begin{equation}
f(z)=z+\sum\limits_{n=2}^{\infty }a_{n}z^{n}  \label{Int-e1}
\end{equation}%
which are analytic in the open unit disc $\mathbb{U}=\{z:z\in \mathbb{C}\,\,%
\mathrm{and}\,\,|z|<1\}.$ Further, by $\mathcal{S}$ we will show the family
of all functions in $\mathcal{A}$ which are univalent in $\mathbb{U}.$

For two functions $f$ and $g,$ analytic in $\mathbb{U},$ we say that the
function $f(z)$ is subordinate to $g(z)$ in $\mathbb{U},$ and write 
\begin{equation*}
f(z) \prec g(z) \qquad \left(z\in \mathbb{U}\right)
\end{equation*}
if there exists a Schwarz function $w(z),$ analytic in $\mathbb{U},$ with 
\begin{equation*}
w(0)=0 \quad \mathrm{and} \quad |w(z)|<1 \qquad \left(z\in \mathbb{U}\right)
\end{equation*}
such that 
\begin{equation*}
f(z)=g(w(z)) \qquad \left(z\in\mathbb{U}\right).
\end{equation*}
In particular, if the function $g$ is univalent in $\mathbb{U},$ the above
subordination is equivalent to 
\begin{equation*}
f(0)=g(0) \quad \mathrm{and} \quad f(\mathbb{U}) \subset g(\mathbb{U}).
\end{equation*}
%%=======Bi-univalent-definition=========

It is well known that every function $f\in \mathcal{S}$ has an inverse $%
f^{-1},$ defined by 
\begin{equation*}
f^{-1}(f(z))=z \qquad \left(z \in \mathbb{U}\right)
\end{equation*}
and 
\begin{equation*}
f(f^{-1}(w))=w \qquad \left(|w| < r_0(f);\,\, r_0(f) \geq \frac{1}{4}\right),
\end{equation*}
where 
\begin{equation}  \label{Int-f-inver}
f^{-1}(w) = w - a_2w^2 + (2a_2^2-a_3)w^3 - (5a_2^3-5a_2a_3+a_4)w^4+\ldots .
\end{equation}

A function $f \in \mathcal{A}$ is said to be bi-univalent in $\mathbb{U}$ if
both $f(z)$ and $f^{-1}(z)$ are univalent in $\mathbb{U}.$ Let $\Sigma$
denote the class of bi-univalent functions in $\mathbb{U}$ given by (\ref%
{Int-e1}). For a brief history and interesting examples of functions which
are in (or which are not in) the class $\Sigma $, together with various
other properties of the bi-univalent function class $\Sigma $ one can refer
the work of Srivastava et al. \cite{HMS-AKM-PG} and references therein. In
fact, the study of the coefficient problems involving bi-univalent functions
was reviewed recently by Srivastava et al. \cite{HMS-AKM-PG}. Various
subclasses of the bi-univalent function class $\Sigma $ were introduced and
non-sharp estimates on the first two coefficients $|a_{2}|$ and $|a_{3}|$ in
the Taylor-Maclaurin series expansion (\ref{Int-e1}) were found in several
recent investigations (see, for example, \cite%
{Ali-Ravi-Ma-Mina-class,Bulut,Caglar-Orhan,Deniz,BAF-MKA,haya,NM-Yamini,GMS-NM-PMK-AAA,SSS-VR-VR,HMS,HMS-Caglar,HMS-GMS-NM-GJM,Xu-HMS-AML,Xu-HMS-AMC}%
). The aforecited all these papers on the subject were actually motivated by
the pioneering work of Srivastava et al. \cite{HMS-AKM-PG}. However, the
problem to find the coefficient bounds on $|a_{n}|$ ($n=3,4,\dots $) for
functions $f\in \Sigma $ is still an open problem.

%%=============================================================

Some of the important and well-investigated subclasses of the univalent
function class $\mathcal{S}$ include (for example) the class $\mathcal{S}%
^{\ast }(\alpha )$ of starlike functions of order $\alpha $ in $\mathbb{U}$
and the class $\mathcal{K}(\alpha )$ of convex functions of order $\alpha $
in $\mathbb{U}.$ By definition, we have 
\begin{equation}
\mathcal{S}^{\ast }(\alpha ):=\left\{ f:f\in \mathcal{A}\,\text{\ }\mathrm{%
and}\,\Re \left( \frac{zf^{\prime }(z)}{f(z)}\right) >\alpha ;\,\,z\in 
\mathbb{U};\,\,\,0\leq \alpha <1\right\}  \label{ST-e}
\end{equation}%
and 
\begin{equation}
\mathcal{K}(\alpha ):=\left\{ f:f\in \mathcal{A}\,\text{\ }\mathrm{and}\,\Re
\left( 1+\frac{zf^{\prime \prime }(z)}{f^{\prime }(z)}\right) >\alpha
;\,\,z\in \mathbb{U};\,\,\,0\leq \alpha <1\right\} .  \label{CV-e}
\end{equation}

For $0\leq \alpha <1,$ a function $f\in\Sigma$ is in the class $%
S^*_{\Sigma}(\alpha)$ of bi-starlike function of order $\alpha,$ or $%
\mathcal{K}_{\Sigma, \alpha}$ of bi-convex function of order $\alpha$ if
both $f$ and $f^{-1}$ are respectively starlike or convex functions of order 
$\alpha.$ For $0< \beta \leq 1,$ a function $f\in\Sigma$ is strongly
bi-starlike function of order $\beta,$ if both the functions $f$ and $f^{-1}$
are strongly starlike of order $\beta.$ We denote the class of all such
functions is denoted by $S^*_{\Sigma, \beta}.$ 
%%====Ma-Minda-starlike and convex=======

Let $\varphi$ be an analytic and univalent function with positive real part
in $\mathbb{U}$ with $\varphi(0)=1,$ $\varphi^{\prime }(0) >0$ and $\varphi$
maps the unit disk $\mathbb{U}$ onto a region starlike with respect to $1,$
and symmetric with respect to the real axis. The Taylor's series expansion
of such function is of the form 
\begin{equation}  \label{varphi-expression}
\varphi(z)=1+B_1z+B_2z^2+B_3z^3+\dots ,
\end{equation}
where all coefficients are real and $B_1 > 0.$ Throughout this paper we
assume that the function $\varphi$ satisfies the above conditions one or
otherwise stated.

By $\mathcal{S}^{\ast }(\varphi )$ and $\mathcal{K}(\varphi )$ we denote the
following classes of functions 
\begin{equation}
\mathcal{S}^{\ast }(\varphi ):=\left\{ f:f\in \mathcal{A}\,\text{\ }\mathrm{%
and}\,\,\frac{zf^{\prime }(z)}{f(z)}\prec \varphi (z);\,\,z\in \mathbb{U}%
\right\}   \label{ST-varphi-e}
\end{equation}%
and 
\begin{equation}
\mathcal{K}(\varphi ):=\left\{ f:f\in \mathcal{A}\text{ }\,\mathrm{and}\,\,1+%
\frac{zf^{\prime \prime }(z)}{f^{\prime }(z)}\prec \varphi (z);\,\,z\in 
\mathbb{U}\right\} .  \label{CV-varphi-e}
\end{equation}

The classes $\mathcal{S}^*(\varphi)$ and $\mathcal{K}(\varphi)$ are the
extensions of a classical sets of a starlike and convex functions and in a
such form were defined and studied by Ma and Minda \cite{Ma-Minda}. A
function $f$ is bi-starlike of Ma-Minda type or bi-convex of Ma-Minda type
if both $f$ and $f^{-1}$ are respectively Ma-Minda starlike or convex. These
classes are denoted respectively by $\mathcal{S}^*_{\Sigma}(\varphi)$ and $%
\mathcal{K}_{\Sigma}(\varphi)$ (see \cite{Ali-Ravi-Ma-Mina-class}). 
%%=======================================================

In order to derive our main results, we will need the following lemma.

\begin{lemma}
\textrm{(see \cite{Pom})}\label{lem-pom} If $p\in \mathcal{P},$ then $%
|p_i|\leq 2$ for each $i,$ where $\mathcal{P}$ is the family of all
functions $p,$ analytic in $\mathbb{U},$ for which 
\begin{equation*}
\Re\{p(z)\}>0\quad (z \in \mathbb{U}),
\end{equation*}
where 
\begin{equation*}
p(z)=1+p_1z+p_2z^2+\cdots \quad (z \in \mathbb{U}).
\end{equation*}
\end{lemma}

%%=================================
Motivated by the aforementioned works (especially \cite{Zaprawa} and \cite%
{Caglar-Orhan,HO-NM-VKB,HMS-Caglar}), we consider the following subclass of
the function class $\Sigma$ (see also, \cite{Tang}). 
%-------main-results---Orhan-class--

A function $f\in \Sigma$ given by (\ref{Int-e1}) is said to be in the class $%
\mathcal{N}^{\mu, \lambda}_{\Sigma}(\varphi)$ if the following conditions
are satisfied: 
\begin{equation}  \label{Orhan-class-e1}
(1-\lambda)\left(\frac{f(z)}{z}\right)^{\mu}+\lambda f^{\prime }(z)\left(%
\frac{f(z)}{z}\right)^{\mu-1} \prec \varphi(z) \qquad (\lambda \geq 1,\,
\mu\geq 0, \, z \in \mathbb{U})
\end{equation}
and 
\begin{equation}  \label{Orhan-class-e2}
(1-\lambda)\left(\frac{g(w)}{w}\right)^{\mu}+\lambda g^{\prime }(w)\left(%
\frac{g(w)}{w}\right)^{\mu-1}\prec \varphi(w) \qquad (\lambda \geq 1,\,
\mu\geq 0, \, w \in \mathbb{U}),
\end{equation}
where $g(w)=f^{-1}(w).$

\begin{remark}
\label{Def-Remark}

From among the many choices of $\mu$, $\lambda$ and the function $\varphi$
which would provide the following known subclasses:

\begin{enumerate}
\item $\mathcal{N}^{1,1}_{\Sigma}(\varphi)$ = $\mathcal{H}%
_{\Sigma}^{\varphi} $ \, \, \cite[p.345]{Ali-Ravi-Ma-Mina-class}.

\item $\mathcal{N}^{1,1}_{\Sigma}(\left(\frac{1+z}{1-z}\right)^{\beta})$ = $%
\mathcal{H}_{\Sigma}^{\beta}$ $(0<\beta\leq 1)$ \, and $\mathcal{N}%
^{1,1}_{\Sigma}(\frac{1+(1-2\alpha)z}{1-z})$ = $\mathcal{H}%
_{\Sigma}^{\alpha} $ \, $(0\leq \alpha <1)$\, \newline
\cite[Definitions 1 and 2]{HMS-AKM-PG}.

\item $\mathcal{N}^{1,\lambda}_{\Sigma}(\varphi)$ = $\mathcal{R}%
_{\Sigma}(\lambda, \varphi)$ \,\, $(\lambda \geq 0)$ \, \cite[Definition 1.1]%
{SSS-VR-VR}.

\item $\mathcal{N}^{1,\lambda}_{\Sigma}(\left(\frac{1+z}{1-z}%
\right)^{\beta}) $ = $\mathcal{B}_{\Sigma}(\beta,\lambda)$ \, $(\lambda \geq
1; 0<\beta\leq 1) $\, and $\mathcal{N}^{1,\lambda}_{\Sigma}(\frac{%
1+(1-2\alpha)z}{1-z})$ = $\mathcal{B}_{\Sigma}(\alpha,\lambda)$\, $(\lambda
\geq 1; 0\leq \alpha <1)$\, \cite[Definitions 2.1 and 3.1]{BAF-MKA}.

\item $\mathcal{N}^{\mu,1}_{\Sigma}(\varphi)$ = $\mathcal{F}%
^{\mu}_{\Sigma}(\varphi)$ \,\, $(\mu \geq 0)$ \, \cite[Definition 2.1]%
{SSS-VR-VR}.

\item $\mathcal{N}^{0,1}_{\Sigma}(\left(\frac{1+z}{1-z}\right)^{\beta})$ = $%
\mathcal{S}^*_{\Sigma,\beta}$ \, $(0<\beta\leq 1)$\, and $\mathcal{N}%
^{0,1}_{\Sigma}(\frac{1+(1-2\alpha)z}{1-z})$ = $\mathcal{S}%
^*_{\Sigma}(\alpha)$ \, $(0\leq \alpha <1).$

\item $\mathcal{N}^{\mu,\lambda}_{\Sigma}(\left(\frac{1+z}{1-z}%
\right)^{\beta})$ = $\mathcal{N}^{\mu,\lambda}_{\Sigma}(\beta)$ \, $(\lambda
\geq 1; \mu \geq 0; 0<\beta\leq 1)$ \, \cite[Definitions 2.1]{Caglar-Orhan}.

\item[ ] and

\item[ ] $\mathcal{N}^{\mu,\lambda}_{\Sigma}(\frac{1+(1-2\alpha)z}{1-z})$ = $%
\mathcal{N}^{\mu,\lambda}_{\Sigma}(\alpha)$ \, $(\lambda \geq 1; \mu \geq 0;
0\leq \alpha <1)$ \, \cite[Definitions 3.1]{Caglar-Orhan}.
\end{enumerate}
\end{remark}

In this paper we shall obtain the Fekete-Szeg\"{o} inequalities for $%
\mathcal{N}^{\mu, \lambda}_{\Sigma}(\varphi)$ and its special classes. These
inequalities will result in bounds of the third coefficient which are, in
some cases, better than these obtained in \cite%
{Ali-Ravi-Ma-Mina-class,Caglar-Orhan,BAF-MKA,HMS-Caglar,HMS-AKM-PG,Tang}.

\section{Main Results}

%%==============theorem-Orhan-class====================

\begin{theorem}
\label{th-orhan-class} Let $f$ of the form (\ref{Int-e1}) be in $\mathcal{N}%
^{\mu, \lambda}_{\Sigma}(\varphi)$ and $\delta \in \mathbb{R}.$ Then 
\begin{equation}  \label{th-orhan-class-a3-a2}
|a_3-\delta a_2^2| \leq \left \{ 
\begin{array}{cc}
\frac{B_1}{2\lambda+\mu} & ; |\delta -1| \leq \frac{\mu+1}{2}\left |1+\frac{%
2(B_1-B_2)(\lambda+\mu)^2}{B_1^2(2\lambda+\mu)(1+\mu)}\right | \\ 
\frac{2B_1^3|\delta - 1|}{|(2\lambda+\mu)(1+\mu)B_1^2+2(B_1-B_2)(\lambda+%
\mu)^2|} & ; |\delta -1| \geq \frac{\mu+1}{2}\left |1+\frac{%
2(B_1-B_2)(\lambda+\mu)^2}{B_1^2(2\lambda+\mu)(1+\mu)}\right |.%
\end{array}%
\right.
\end{equation}
\end{theorem}

\begin{proof}
Since $f\in \mathcal{N}^{\mu, \lambda}_{\Sigma}(\varphi),$ there exists two
analytic functions $r,s:\mathbb{U} \rightarrow \mathbb{U},$ with $%
r(0)=0=s(0),$ such that 
\begin{equation}  \label{th-orhan-class-p-e3}
(1-\lambda)\left(\frac{f(z)}{z}\right)^{\mu}+\lambda f^{\prime }(z)\left(%
\frac{f(z)}{z}\right)^{\mu-1}=\varphi(r(z))
\end{equation}
and 
\begin{equation}  \label{th-orhan-class-p-e4}
(1-\lambda)\left(\frac{g(w)}{w}\right)^{\mu}+\lambda g^{\prime }(w)\left(%
\frac{g(w)}{w}\right)^{\mu-1}=\varphi(s(z)).
\end{equation}
Define the functions $p$ and $q$ by 
\begin{equation}  \label{th-orhan-class-p-e5}
p(z)=\frac{1+r(z)}{1-r(z)}=1+p_1z+p_2z^2+p_3z^3+\dots
\end{equation}
and 
\begin{equation}  \label{th-orhan-class-p-e6}
q(z)=\frac{1+s(z)}{1-s(z)}=1+q_1z+q_2z^2+q_3z^3+\dots
\end{equation}
or equivalently, 
\begin{equation}  \label{th-orhan-class-p-e7}
r(z)=\frac{p(z)-1}{p(z)+1}=\frac{1}{2}\left(p_1z+\left(p_2-\frac{p_1^2}{2}%
\right)z^2+\left(p_3+\frac{p_1}{2}\left(\frac{p_1^2}{2}-p_2\right)-\frac{%
p_1p_2}{2}\right)z^3+\dots\right)
\end{equation}
and 
\begin{equation}  \label{th-orhan-class-p-e8}
s(z)=\frac{q(z)-1}{q(z)+1}=\frac{1}{2}\left(q_1z+\left(q_2-\frac{q_1^2}{2}%
\right)z^2+\left(q_3+\frac{q_1}{2}\left(\frac{q_1^2}{2}-q_2\right)-\frac{%
q_1q_2}{2}\right)z^3+\dots\right).
\end{equation}

Using (\ref{th-orhan-class-p-e7}) and (\ref{th-orhan-class-p-e8}) in (\ref%
{th-orhan-class-p-e3}) and (\ref{th-orhan-class-p-e4}), we have 
\begin{equation}  \label{th-orhan-class-p-e3-9}
(1-\lambda)\left(\frac{f(z)}{z}\right)^{\mu}+\lambda f^{\prime }(z)\left(%
\frac{f(z)}{z}\right)^{\mu-1}=\varphi\left(\frac{p(z)-1}{p(z)+1}\right)
\end{equation}
and 
\begin{equation}  \label{th-orhan-class-p-e3-10}
(1-\lambda)\left(\frac{g(w)}{w}\right)^{\mu}+\lambda g^{\prime }(w)\left(%
\frac{g(w)}{w}\right)^{\mu-1}=\varphi\left(\frac{q(w)-1}{q(w)+1}\right).
\end{equation}

Again using (\ref{th-orhan-class-p-e7}) and (\ref{th-orhan-class-p-e8})
along with (\ref{varphi-expression}), it is evident that

\begin{equation}  \label{th-orhan-class-p-e11}
\varphi\left(\frac{p(z)-1}{p(z)+1}\right)=1+\frac{1}{2}B_1p_1z+\left(\frac{1%
}{2}B_1\left(p_2-\frac{1}{2}p_1^2\right)+\frac{1}{4}B_2p_1^2\right)z^2 +
\dots
\end{equation}
and 
\begin{equation}  \label{th-orhan-class-p-e12}
\varphi\left(\frac{q(w)-1}{q(w)+1}\right) =1+\frac{1}{2}B_1q_1w+\left(\frac{1%
}{2}B_1\left(q_2-\frac{1}{2}q_1^2\right)+\frac{1}{4}B_2q_1^2\right)w^2 +
\dots.
\end{equation}

It follows from (\ref{th-orhan-class-p-e3-9}), (\ref{th-orhan-class-p-e3-10}%
), (\ref{th-orhan-class-p-e11}) and (\ref{th-orhan-class-p-e12}) that 
\begin{equation}  \label{th-orhan-class-p-e13}
(\lambda+\mu)a_2 = \frac{1}{2}B_1p_1
\end{equation}
\begin{equation}  \label{th-orhan-class-p-e14}
(2\lambda+\mu)[a_3+\frac{a_2^2}{2}(\mu-1)]= \frac{1}{2}B_1\left(p_2-\frac{1}{%
2}p_1^2\right)+\frac{1}{4}B_2p_1^2
\end{equation}
\begin{equation}  \label{th-orhan-class-p-e15}
-(\lambda+\mu)a_2 = \frac{1}{2}B_1q_1
\end{equation}
and 
\begin{equation}  \label{th-orhan-class-p-e16}
(2\lambda+\mu)[\frac{a_2^2}{2}(\mu+3)-a_3]= \frac{1}{2}B_1\left(q_2-\frac{1}{%
2}q_1^2\right)+\frac{1}{4}B_2q_1^2.
\end{equation}
From (\ref{th-orhan-class-p-e13}) and (\ref{th-orhan-class-p-e15}), we find
that 
\begin{equation}  \label{orhan_a-2}
a_2 = \frac{B_1p_1}{2(\lambda+\mu)}=\frac{-B_1q_1}{2(\lambda+\mu)}
\end{equation}
it follows that 
\begin{equation}  \label{th-orhan-class-p-e17}
p_1 =-q_1
\end{equation}
and 
\begin{equation}  \label{th-orhan-class-p-e18}
8(\lambda+\mu)^2a_2^2=B_1^2(p_1^2+q_1^2).
\end{equation}
Adding (\ref{th-orhan-class-p-e14}) and (\ref{th-orhan-class-p-e16}), we
have 
\begin{equation}  \label{orhan-a-2-square}
a_2^2(2\lambda+\mu)(\mu+1)=\frac{B_1}{2}(p_2+q_2)+\frac{(B_2-B_1)}{4}%
(p_1^2+q_1^2).
\end{equation}
Substituting (\ref{orhan_a-2}) and (\ref{th-orhan-class-p-e17}) into (\ref%
{orhan-a-2-square}), we get, 
\begin{equation}  \label{p-2-square}
p_1^2=\frac{B_12(\lambda+\mu)^2(p_2+q_2)}{B_1^2(2\lambda+\mu)(%
\mu+1)-2(B_2-B_1)(\lambda+\mu)^2}.
\end{equation}
Now, (\ref{orhan_a-2}) and (\ref{p-2-square}) yield 
\begin{equation}  \label{th-orhan-class-Fekete-a2}
a_2^2=\frac{B_1^3(p_2+q_2)}{2(\mu+1)(2\lambda+\mu)B_1^2+4(B_1-B_2)(\lambda+%
\mu)^2}.
\end{equation}

By subtracting (\ref{th-orhan-class-p-e14}) from (\ref{th-orhan-class-p-e16}%
) and a computation using (\ref{th-orhan-class-p-e17}) finally lead to 
\begin{equation}
a_{3}=a_{2}^{2}+\frac{B_{1}(p_{2}-q_{2})}{8\lambda +4\mu }.
\label{th-orhan-class-Fekete-a3}
\end{equation}%
From (\ref{th-orhan-class-Fekete-a2}) and (\ref{th-orhan-class-Fekete-a3})
it follows that 
\begin{equation*}
a_{3}-\delta a_{2}^{2}=B_{1}\left[ \left( h(\delta )+\frac{1}{8\lambda +4\mu 
}\right) p_{2}+\left( h(\delta )-\frac{1}{8\lambda +4\mu }\right) q_{2}%
\right] ,
\end{equation*}%
where 
\begin{equation*}
h(\delta )=\frac{B_{1}^{2}(1-\delta )}{2(\mu +1)(2\lambda +\mu
)B_{1}^{2}+4(B_{1}-B_{2})(\lambda +\mu )^{2}}.
\end{equation*}%
Since all $B_{j}$ are real and $B_{1}>0,$ we conclude that 
\begin{equation*}
|a_{3}-\delta a_{2}^{2}|\leq \left\{ 
\begin{array}{cc}
\frac{B_{1}}{2\lambda +\mu } & ;0\leq |h(\delta )|<\frac{1}{8\lambda +\mu }
\\ 
4B_{1}|h(\delta )| & ;|h(\delta )|\geq \frac{1}{8\lambda +\mu }~,%
\end{array}%
\right.
\end{equation*}%
which completes the proof.
\end{proof}

\begin{remark}
For $\lambda=\mu=1$ Theorem \ref{th-orhan-class} reduces to the results
discussed in \cite[Theorem 1, p.172]{Zaprawa}.
\end{remark}

%%-----------------Special-cases------------------------------

\section{Corollaries and Consequences}

Taking $\delta=1,$ $\delta=0$ in Theorem \ref{th-orhan-class}, we have the
following corollaries.

\begin{corollary}
\label{fekete-paper-Cor1} If $f\in \mathcal{N}^{\mu,
\lambda}_{\Sigma}(\varphi)$ then 
\begin{equation*}
|a_3-a_2^2| \leq \frac{B_1}{2\lambda+\mu}.
\end{equation*}
\end{corollary}

\begin{corollary}
\label{fekete-paper-Cor2} If $f\in \mathcal{N}^{\mu,
\lambda}_{\Sigma}(\varphi)$ then 
\begin{equation*}
|a_3| \leq \left \{ 
\begin{array}{cc}
\frac{B_1}{2\lambda+\mu} & ; \tfrac{(B_1-B_2)}{B_1^2} \in \left (-\infty, 
\tfrac{-(3+\mu)(2\lambda+\mu)}{2(\lambda+\mu)^2}\right]\bigcup \left[\tfrac{%
(1-\mu)(2\lambda+\mu)}{2(\lambda+\mu)^2}, \infty \right ) \\ 
\tfrac{2B_1^3}{|(2\lambda+\mu)(1+\mu)B_1^2+2(B_1-B_2)(\lambda+\mu)^2|} & ; 
\tfrac{(B_1-B_2)}{B_1^2} \in \left [\tfrac{-(3+\mu)(2\lambda+\mu)}{%
2(\lambda+\mu)^2}, \tfrac{-(1+\mu)(2\lambda+\mu)}{2(\lambda+\mu)^2}%
\right)\bigcup \left(\tfrac{-(1+\mu)(2\lambda+\mu)}{2(\lambda+\mu)^2}, 
\tfrac{(1-\mu)(2\lambda+\mu)}{2(\lambda+\mu)^2}\right ].%
\end{array}%
\right.
\end{equation*}
\end{corollary}

\begin{remark}
Corollary \ref{fekete-paper-Cor2} provides an improvement of the estimate $%
|a_3|$ obtained by Tang et al. \cite[Theorem 2.1, p.3]{Tang}.
\end{remark}

In view of Remark \ref{Def-Remark}, Corollaries \ref{fekete-paper-Cor1} and %
\ref{fekete-paper-Cor2} yield the following corollaries.

\begin{corollary}
\label{fekete-paper-Cor3} If $f\in \mathcal{N}^{\mu,
\lambda}_{\Sigma}(\beta) $ then 
\begin{equation*}
|a_3| \leq \frac{2\beta}{2\lambda+\mu} \qquad \mathrm{and} \qquad
|a_3-a_2^2| \leq \frac{2\beta}{2\lambda+\mu}.
\end{equation*}
\end{corollary}

\begin{corollary}
\label{fekete-paper-Cor4} If $f\in \mathcal{N}^{\mu,
\lambda}_{\Sigma}(\alpha)$ then 
\begin{equation*}
|a_3| \leq \frac{2(1-\alpha)}{2\lambda+\mu} \qquad \mathrm{and} \qquad
|a_3-a_2^2| \leq \frac{2(1-\alpha)}{2\lambda+\mu}.
\end{equation*}
\end{corollary}

\begin{remark}
The bounds $|a_3|$ obtained in Corollaries \ref{fekete-paper-Cor3} and \ref%
{fekete-paper-Cor4} are improvement of the bounds $|a_3|$ estimated by \c{C}a%
\u{g}lar et al. \cite[Theorems 2.1 and 3.1]{Caglar-Orhan}.
\end{remark}

%%===========================

\begin{remark}
In view of Remark \ref{Def-Remark} the aforecited work for the subclasses $%
\mathcal{H}^{\varphi}_{\Sigma},$ $\mathcal{H}^{\beta}_{\Sigma}$ and $%
\mathcal{H}^{\alpha}_{\Sigma}$ are coincide with the results of Zaprawa \cite%
[Corollaries 1 to 4, p.173]{Zaprawa}.
\end{remark}

%%===============Ma-Minda-Starlike-Cases=============================

\begin{corollary}
\label{fekete-paper-Cor11} If $f\in \mathcal{S}^*_{\Sigma}(\varphi)$ then 
\begin{equation*}
|a_3-a_2^2| \leq \frac{B_1}{2}.
\end{equation*}
\end{corollary}

\begin{corollary}
\label{fekete-paper-Cor12} If $f\in \mathcal{S}^*_{\Sigma}(\varphi)$ then 
\begin{equation*}
|a_3| \leq \left \{ 
\begin{array}{cc}
\frac{B_1}{2} & ; \frac{(B_1-B_2)}{B_1^2} \in \left (-\infty, -3\right]%
\bigcup \left[0, \infty \right) \\ 
\frac{B_1^3}{|B_1^2+(B_1-B_2)|} & ; \frac{(B_1-B_2)}{B_1^2} \in \left [-2,
-1\right)\bigcup \left(-1, 1\right ].%
\end{array}%
\right.
\end{equation*}
\end{corollary}

\begin{corollary}
\label{fekete-paper-Cor13} If $f\in \mathcal{S}^*_{\Sigma, \beta}$ then 
\begin{equation*}
|a_3| \leq \beta \qquad \mathrm{and } \qquad |a_3-a_2^2| \leq \beta .
\end{equation*}
\end{corollary}

\begin{corollary}
\label{fekete-paper-Cor14} If $f\in \mathcal{S}^*_{\Sigma}(\alpha)$ then 
\begin{equation*}
|a_3| \leq 1-\alpha \qquad \mathrm{and } \qquad |a_3-a_2^2| \leq 1-\alpha .
\end{equation*}
\end{corollary}

\begin{remark}
The inequalities estimated in Corollaries \ref{fekete-paper-Cor12} to \ref%
{fekete-paper-Cor14} are improvement of the inequalities obtained by Zaprawa 
\cite[Corollaries 11 and 12, p.174]{Zaprawa}.
\end{remark}

%%=====================Aouf-Frasin-Class==========================

\begin{corollary}
\label{fekete-paper-Cor16} If $f\in \mathcal{R}_{\Sigma}(\lambda;\varphi)$
then 
\begin{equation*}
|a_3-a_2^2| \leq \frac{B_1}{2\lambda+1}.
\end{equation*}
\end{corollary}

\begin{corollary}
\label{fekete-paper-Cor17} If $f\in \mathcal{R}_{\Sigma}(\lambda;\varphi)$
then 
\begin{equation*}
|a_3| \leq \left \{ 
\begin{array}{cc}
\frac{B_1}{2\lambda+1} & ; \frac{(B_1-B_2)}{B_1^2} \in \left (-\infty, \frac{%
2(2\lambda+1)}{(\lambda+1)^2}\right]\bigcup \left[0, \infty \right) \\ 
\frac{B_1^3}{|(2\lambda+1)B_1^2+(B_1-B_2)(\lambda+1)^2|} & ; \frac{(B_1-B_2)%
}{B_1^2} \in \left [\frac{2(2\lambda+1)}{(\lambda+1)^2}, \frac{-(2\lambda+1)%
}{(\lambda+1)^2}\right)\bigcup \left(\frac{-(2\lambda+1)}{(\lambda+1)^2},
0\right ] .%
\end{array}%
\right.
\end{equation*}
\end{corollary}

\begin{remark}
Corollary \ref{fekete-paper-Cor17} provides an improvement of $|a_3|$
obtained by Sivaprasad Kumar et al. \cite[Theorem 2.1, p.3]{SSS-VR-VR}.
\end{remark}

\begin{corollary}
\label{fekete-paper-Cor18} If $f\in \mathcal{B}_{\Sigma}(\beta, \lambda)$
then 
\begin{equation*}
|a_3| \leq \frac{2\beta}{2\lambda+1} \qquad \mathrm{and } \qquad |a_3-a_2^2|
\leq \frac{2\beta}{2\lambda+1}.
\end{equation*}
\end{corollary}

\begin{corollary}
\label{fekete-paper-Cor19} If $f\in \mathcal{B}_{\Sigma}(\alpha, \lambda)$
then 
\begin{equation*}
|a_3| \leq \frac{2(1-\alpha)}{2\lambda+1} \qquad \mathrm{and } \qquad
|a_3-a_2^2| \leq \frac{2(1-\alpha)}{2\lambda+1}.
\end{equation*}
\end{corollary}

\begin{remark}
The bounds $|a_{3}|$ obtained in Corollaries \ref{fekete-paper-Cor18} and %
\ref{fekete-paper-Cor19} are improvement of the bounds $|a_{3}|$ estimated
by Frasin and Aouf \cite[Theorems 2.2 and 3.2, p.1570 and 1572]{BAF-MKA},
respectively.
\end{remark}

%%===========================================================

\begin{remark}
\label{varphi-0} If we take 
\begin{equation}  \label{e-varphi-0}
\varphi=\varphi_0=\frac{1+z}{1-z}=1+2z+2z^2+\dots
\end{equation}
in the class $\mathcal{N}^{\mu, \lambda}_{\Sigma}(\varphi),$ we are led to
the class which we denote, for convenience, by $\mathcal{N}^{\mu,
\lambda}_{\Sigma}(\varphi_0).$ In particular, $\mathcal{N}^{1,
1}_{\Sigma}(\varphi_0)=:\mathcal{H}^{\varphi_0}_{\Sigma},$ $\mathcal{N}^{0,
\mu}_{\Sigma}(\varphi_0)=:\mathcal{S}^*_{\Sigma}(\varphi_0)$ and $\mathcal{N}%
^{1, \lambda}_{\Sigma}(\varphi)=:\mathcal{B}^*_{\Sigma}(\lambda,\varphi_0).$
\end{remark}

In view of Remark \ref{varphi-0}, the Corollaries \ref{fekete-paper-Cor1}
and \ref{fekete-paper-Cor2} yield the following corollaries.

\begin{corollary}
\label{fekete-paper-Cor5} If $f\in \mathcal{N}^{\mu,
\lambda}_{\Sigma}(\varphi_0)$ then 
\begin{equation*}
|a_3| \leq \frac{2}{2\lambda+\mu} \qquad \mathrm{and} \qquad |a_3-a_2^2|
\leq \frac{2}{2\lambda+\mu}.
\end{equation*}
\end{corollary}

\begin{remark}
For $\mu=\lambda=1$ the estimates in Corollary \ref{fekete-paper-Cor5} would
reduce to a known result in \cite[Corollary 5, p.173]{Zaprawa}
\end{remark}

\begin{corollary}
\label{fekete-paper-Cor15} If $f\in \mathcal{S}^*_{\Sigma}(\varphi_0)$ then 
\begin{equation*}
|a_3| \leq 1 \qquad \mathrm{and } \qquad |a_3-a_2^2| \leq 1.
\end{equation*}
\end{corollary}

\begin{corollary}
\label{fekete-paper-Cor20} If $f\in \mathcal{B}_{\Sigma}(\lambda, \varphi_0)$
then 
\begin{equation*}
|a_3| \leq \frac{2}{2\lambda+1} \qquad \mathrm{and } \qquad |a_3-a_2^2| \leq 
\frac{2}{2\lambda+1}.
\end{equation*}
\end{corollary}

%-------References------------

\end{document}